\documentclass[12pt]{article} 
\pdfoutput=1
\usepackage[a4paper]{geometry}
\usepackage[utf8]{inputenc} \usepackage[english]{babel} \usepackage[T2A]{fontenc}
\usepackage{amsmath,amssymb} \usepackage{euler,eucal}
\usepackage{myconcrete} 
\usepackage[unicode=true,plainpages=false]{hyperref} 
\hypersetup{colorlinks=true,linkcolor=magenta,anchorcolor=magenta,urlcolor=blue,citecolor=blue}
\usepackage{algorithmic} 
\usepackage[ruled]{algorithm} 
\usepackage{graphicx}
\usepackage{setspace}
\def\title#1{\phantom{.}\vskip 5mm
      \begin{center}\begin{doublespace}{ \LARGE\sc  #1}\end{doublespace}\end{center}}

\def\author#1{\begin{center}\large{#1}\end{center}}

\def\address#1{\begin{center}\footnotesize\textit{#1}\end{center}}

\def\date#1{\vskip 3mm\begin{center}\large{#1}\end{center}\vskip 10mm}

\def\affilnum#1{${}^{\@fnsymbol{#1}}$}
\def\affil#1{${}^{\@fnsymbol{#1}}$}

\def\thanks#1{${}^{\@crssymbol{#1}}$}
\def\grant#1{\def\thefootnote{\$}\footnotetext{#1} }

\def\email#1{\texttt{#1}}
\def\keywordsname{Keywords}
\def\aabstractname{Abstract.}
\renewenvironment{abstract}{\small \quotation {\bfseries \aabstractname} }{\endquotation}
\def\keywords#1{\vskip 2mm\par\noindent{\small\emph{\keywordsname:} #1}}
\usepackage{caption}
\usepackage{amsthm}

\theoremstyle{definition}
\newtheorem{remark}{Remark}
\newtheorem{corollary}{Corollary}
\newtheorem{theorem}{Theorem}
\def\hat{\mathaccent"705E }
\let\eps=\varepsilon 
  
\let\leq=\leqslant 
\def\pptitle{Fast truncation of mode ranks for bilinear tensor operations}
\def\ppauthor{D. V. Savostyanov, E. E. Tyrtyshnikov, N. L. Zamarashkin}

\def\ppkeywords{%
          Multidimensional arrays, structured tensors, Tucker approximation, 
          fast compression, cross approximation
          }

\hypersetup{%
            pdfauthor={\ppauthor},
            pdftitle={\pptitle},
            pdfsubject={Research article on tensor approximations},
            pdfkeywords={\ppkeywords}
            }

\def\trans{T} \def\t{\trans}

\def\eps{\varepsilon} \def\phi{\varphi}

\def\A{\mathbf{A}} \def\B{\mathbf{B}} 
\def\U{\mathbf{U}}  
\def\F{\mathbf{F}} \def\G{\mathbf{G}} \def\H{\mathbf{H}}
\def\T{\mathbf{T}}
\def\u{\mathbf{u}} \def\v{\mathbf{v}} \def\w{\mathbf{w}}
\def\a{\mathbf{a}} \def\b{\mathbf{b}}  \def\d{\mathbf{d}}
\def\xx{\mathbf{x}} 

\def\O{\mathcal{O}}

\def\x{\mathbin{\times}}
\def\o{\mathbin{\raise1pt\hbox{$\scriptscriptstyle\mathord\otimes$}}}
\def\hadamar{\mathbin{\raise1pt\hbox{$\scriptstyle\mathord\odot$}}}

\newcommand{\eee}[1]{_{\cdot 10}{#1}}

\def\nrm{\mathop{\mathtt{nrm}}\nolimits}
\def\err{\mathop{\mathtt{err}}\nolimits}
\def\rmax{r_\mathrm{max}}
\def\tol{\mathtt{tol}}

\def\rrr{{r_1 \times r_2 \times r_3}}
\def\ppp{{p_1 \times p_2 \times p_3}}
\def\nnn{{n_1 \times n_2 \times n_3}}

\renewcommand{\f}[2]{\left\langle {#1},\:{#2} \right\rangle}
\def\Kron{\mathop{\mathbf{Kron}}\nolimits}

\def\Span{\mathop{\mathrm{span}}\nolimits}

\def\Sum{\sum\limits}
\def\eqdef{\mathrel{\stackrel{\mathrm{def}}{=}}}
\def\tn{\mathbf{2}}

\begin{document}
\bibliographystyle{hsiam}

\title{\pptitle}
\author{\ppauthor}
\address{Institute of Numerical Mathematics, Russian Academy of Sciences, \\
                     Russia, 119333 Moscow, Gubkina 8 \\ 
                  \email{dmitry.savostyanov@gmail.com, [tee,kolya]@bach.inm.ras.ru}
           }
\grant{This work was supported by 
         RFBR grants 08-01-00115, 09-01-12058, 10-01-00757, 10-01-09201, RFBR/DFG grant 09-01-91332,
         Russian Federation Gov. contract $\Pi940$
         and  Priority research program of Dep. Math. RAS.
         }
\date{\today}

\begin{abstract}
We propose a fast algorithm for mode rank truncation of the result of a bilinear operation on 3-tensors given in the Tucker or canonical form. 
If the arguments and the result have mode sizes $n$ and mode ranks $r,$ the computation costs~$\O(nr^3 + r^4).$ 
The algorithm is based on the cross approximation of Gram matrices, and the accuracy of the resulted Tucker approximation is limited by square root of machine precision.
\end{abstract}

\section{Introduction}
Data sparse representations of tensors and efficient operations in the corresponding formats play increasingly important role in many applications. 
In the paper we consider a 3-\emph{tensor} $\A = \A[i,j,k]$ that is an array with three indices. 
The number of allowed values of each index is called \emph{mode size}.
To specify tensor indices explicitly, we use \emph{square brackets}. 
This notation allows to easily specify different \emph{index transformations}.
For instance, \emph{unfoldings} of $\nnn$ tensor $\A[i,j,k],$ are \emph{matricizations} of sizes $n_1 \x n_2n_3,$ $n_2 \x n_1n_3 $ and $n_3 \x n_1n_2$ that consist of columns, rows and tube fibers of $\A,$ 
\begin{equation}\label{eq:un}
A^{(1)} = A[i, jk], \qquad A^{(2)} = A[j, ki], \qquad A^{(3)} = A[k, ij].
\end{equation}
Here we set row/column/fiber index of the tensor $\A[i,j,k]$ as row index and join the two others in one multiindex for columns of the unfolding. 
The result is considered as a two-index object (matrix), with row and column indices separated by comma. 
The difference between matrices and tensors is additionally stressed by use of uppercase letter instead of bold uppercase.
The \emph{reshape} of tensor elements assumes as well a change of the index ordering. 
For example, transposition of matrix reads $(A[i,j])^\t = A[j,i],$ vectorization reads $\a[ij] = A[i,j].$ 
We see that the square bracket notation is rather self-explaining and suits for description of algorithms working with multidimensional data. 
We also will use the MATLAB-style \emph{round bracket} notation $a(i,j,k)$ to point to individual element of $\A[i,j,k]$ and $\a(i,:,k)$ to select a mode vector (i.e. row) from tensor $\A.$ 
Scalars and vectors are denoted by lowercase and bold lowercase letters.

In numerical work with tensors of large \emph{mode size} it is crucial to use data sparse formats. 
For 3-tensors, the most useful are the following.

The \emph{canonical decomposition}  \cite{hitchcock-sum-1927,cc-parafac-1970,harshman-parafac-1970} (or \emph{canonical approximation} to some other tensor) reads
\begin{equation}\label{eq:c}
\A[i,j,k] = \Sum_{s=1}^R \u_s[i] \o \v_s[j] \o \w_s[k], \qquad   a(i,j,k) = \Sum_{s=1}^R u(i,s) v(j,s) w(k,s).
\end{equation}
The minimal possible number of summands is called \emph{tensor rank} or \emph{canonical rank} of the given tensor $\A.$
However, canonical decomposition/approximation of a tensor with minimal value of $R$ is a rather ill-posed and computationally unstable problem \cite{desilva-2008}. 
This explains why among many algorithms of canonical approximation (cf.~\cite{comon-2000,lathauwer-schur-2004,esgras-bb-2009,ost-sorto-2009}) none is known as absolutely reliable, and no robust tools for linear algebra operations maintaining the canonical format (linear combinations, etc.)  are proposed.

The (truncated
) \emph{Tucker decomposition/approximation} \cite{Tucker} reads 
\begin{equation}\label{eq:t}
 \begin{split}
      \A[i,j,k] & = \G[p,q,s] \x_1 U[i,p] \x_2 V[j,q] \x_3 W[k,s], \\
       a(i,j,k) & = \Sum_{p=1}^{r_1}\Sum_{q=1}^{r_2} \Sum_{s=1}^{r_3} g(p,q,s) u(i,p) v(j,q) w(k,s).
 \end{split} 
\end{equation}
The quantities $r_1,r_2,r_3$ are referred to as \emph{Tucker ranks} or \emph{mode ranks},
the tensor $\G=\G[p,q,s]$ of size $r_1 \x r_2 \x r_3$ is called the \emph{Tucker core},
the symbol $\x_l$ designates the multiplication of a tensor by a matrix along the $l$-th mode, 
the \emph{mode factors} $U, V, W$ have orthonormal columns.
In $d$ dimensions, the memory to store $r\x r\x \ldots \x r$ core is $r^d,$ that is usually beyond affordable for large $d$ and even for very small $r$  (so-called \emph{curse of dimensionality}). For $d=3, \: r\sim 100$ the storage is small and Tucker decomposition can be used efficiently.

In~\cite{ost-latensor-2009} the efficient operations with 3-tensors in canonical and Tucker formats are discussed, with approximation of the result in the Tucker format.
Simple operations like linear combination of small number of structured tensors can be done using \emph{multilinear SVD} \cite{lathauwer-svd-2000} (or high-order SVD, HOSVD), with quasi-optimal ranks and guaranteed accuracy.
Linear combination of many tensors, convolution, Hadamard (pointwise) product of tensors and many other bilinear operations reduce to recompression of the following structured tensor
\begin{equation}\label{eq:f}
 \begin{split}
  \F[i,j,k] = & \Kron(\G, \H)[ap,bq,cs] \x_1 U[i,ap] \x_2 V[j,bq] \x_3 W[k,cs], \\
  f(i,j,k)  = & \sum_{pqs} \sum_{abc} g(p,q,s) h(a,b,c) u(i,ap) v(j,bq) w(k,cs),
 \end{split}
\end{equation}
with $\rrr$ core $\G[p,q,s],$  $\ppp$ core $\H[a,b,c]$ and non-orthogonal factors $U, V$ and $W.$
Formally~\eqref{eq:f} is a Tucker-like format with larger mode ranks $p_1r_1, p_2r_2, p_3r_3,$ that should be reduced (truncated) maintaining the desired accuracy.
Due to memory limitations, $\F[i,j,k]$  can not be assembled for mode sizes $n \gtrsim 10^3$ and auxiliary $p_1r_1 \x p_2r_2 \x p_3r_3$ core can not be assembled for ranks $r \gtrsim 30$ (see Tab.~\ref{tab:mem}).\footnote{We always assume $n_1=n_2=n_3=n$ and $r_1=r_2=r_3=p_1=p_2=p_3=r$ in complexity estimates}
The structure of $\F$ should be exploited without explicit evaluation of large temporary arrays.

\begin{table}[t]
\caption{Memory for $r^d$ elements, MB} \label{tab:mem}
\begin{center}
\begin{tabular}{c|cccc}
       &  $d=3$     &  $d=4$         & $d=5$      & $d=6$    \\ \hline
$r=15$ &  $0.026$   &  $0.4$         & $5.8$      & $87$     \\
$r=30$ &  $0.2$     &  $6.2$         & $185$      & $5560$   \\
$r=50$ &  $0.95$    &  $47$          & $2384$     & $119210$  \\
$r=100$&  $7.7$     &  $763$         & $76300$    & $\approx 8$ TB  \\
\end{tabular}
\end{center}
\end{table}

A practical rank-reduction algorithm proposed in \cite{ost-latensor-2009} is a rank revealing version of iterative Tucker-ALS~\cite{tuckerals-1980,lathauwer-rank1-2000} requiring $\O(nr^4 + r^6)$ operations. 
However, the number of iterations in Tucker-ALS depends on the initial guess, and fast approximate evaluation of Tucker factors of~\eqref{eq:f} is important. 

In Sec.~2 we propose to approximate dominant mode subspaces of $\F[i,j,k]$ by the ones of simpler tensors.
In Sec.~3 we compute dominant mode subspaces by a cross approximation of Gram matrices of the unfoldings.
The resulted algorithm requires $\O(nr^3 + r^4)$ operations in three-dimensional case and can be easily generalized to higher dimensions using $\O(dnr^3 + dr^{d+1})$ operations. Since it uses decomposition of Gram matrices, the accuracy  is limited by square root of machine precision. 
In Sec.~4 we apply the proposed method to Hadamard product of electron densities of simple molecules and show that using the result as an initial guess, Tucker-ALS converges to almost machine precision in one iteration.

In the paper we use Frobenius norm of tensors, that is defined as follows
$$
\|\A \|_F^2 \eqdef \f{\A}{\A}, \qquad \f{\A}{\B} \eqdef \Sum_{i=1}^{n_1} \Sum_{j=1}^{n_2} \Sum_{k=1}^{n_3} a_{ijk} b_{ijk}
$$
and spectral norm of tensor (cf.~\cite{defant-1993})
$$
\|\A\|_\tn \eqdef \max_{\|\u\|=\|\v\|=\|\w\|=1} \A \x_1 \u^\t \x_2 \v^\t \x_3 \w^\t = \max_{\|\u\|=\|\v\|=\|\w\|=1} \f{\A}{\u \o \v \o \w},
$$
induced by standard vector norm $\|\u\|^2 \eqdef \|\u\|_2^2  =  (\u,\u) = \sum_{i=1}^n |u_i|^2.$

\section{Approximation of dominant subspaces}
Our goal is to find approximate dominant subspaces of an $\nnn$ tensor~\eqref{eq:f} producing an approximation in the Tucker format
\begin{equation}\label{eq:appr}
 \tilde\F[i,j,k] = \T[\alpha,\beta,\gamma] \x_1 X[i,\alpha] \x_2 Y[j,\beta] \x_3 Z[k,\gamma], \qquad \|\F - \tilde \F\|_F \leq \eps \|\F\|_F
\end{equation}
with a desired (not very high) accuracy and values of mode ranks for $\tilde\F,$ close to optimal.

Tucker factors $X[i,\alpha], Y[j,\beta]$ and $Z[k,\gamma]$ approximate dominant subspaces of rows, columns and fibers of $\F[i,j,k],$ respectively. 
They can be computed by SVD of the unfoldings of $\F,$ as proposed in~\cite{lathauwer-svd-2000}, but this method requires evaluation of all elements of tensor and is not feasible for large mode sizes. 
We can compute~\eqref{eq:appr} interpolating a given tensor on carefully selected set of elements. 
This is done in Cross3D algorithm~\cite{ost-tucker-2008}, that requires evaluation of  $\O(nr + r^3)$ tensor elements and uses $\O(nr^2 + r^4)$ additional operations. 
For a structured tensor~\eqref{eq:f} this summarizes to $\O(nr^3 + r^6)$ operations, i.e. the complexity is~\emph{linear} in mode size.
However, pivoting and error checking involves heuristics and in certain cases is slower than the approximation itself. 
For example, computation of residual $(\F - \tilde\F)[i,j,k]$ on $\O(n)$ randomly picked elements uses $\O(nr^4)$ operations.

To avoid heuristic approaches, we can evaluate dominant subspaces by proper decomposition of Gram matrices of the unfoldings. 
In~\cite{sav-rr-2009} this idea was used for fast mode rank truncation of tensor given in the canonical form~\eqref{eq:c} with large number of terms.
The proposed in~\cite{sav-rr-2009} cross approximation algorithm is equivalent to an unfinished Cholesky decomposition and computes rank-$r$ dominant basis using the diagonal and certain $r$ columns of the Gram matrix.
However, for the unfolding $F[i,jk]$ of tensor $\F[i,j,k]$ the Gram matrix $(F F^\t)[i,i']$ reads
\begin{equation}\label{eq:ff}
 \begin{split}
  (F F^\t)(i,i') = \sum_{pqs} \sum_{abc} \sum_{p'q's'} \sum_{a'b'c'} & g(p,q,s) h(a,b,c) g(p',q',s') h(a',b',c') \\
                                                                     & (V^\t V)(bq,b'q') (W^\t W)(cs,c's') u(i,ap) u(i',a'p'),
 \end{split}
\end{equation}
and it is easy to check, that evaluation of any element of~\eqref{eq:ff} requires $\O(r^6)$ operations.
Therefore, the algorithm from~\cite{sav-rr-2009} applied to~\eqref{eq:ff} has $\O(nr^6)$ complexity, which is not promising even for moderate $r.$
To perform faster, we propose to change the computational objective and look for dominant subspaces of tensors with a simpler structure.

Rewrite the tensor~\eqref{eq:f} as follows
\begin{equation}\label{eq:split}
 \begin{split}
 \F[i,j,k] = & {} \U'[i,bq,cs] \x_2 V[j,bq] \x_3 W[k,cs], \\
             & {} \U'[i,bq,cs] = \Kron(\G,\H)[ap,bq,cs] \x_1 U(i,ap).
 \end{split}
\end{equation}
It is clear that the Tucker approximation of $\U'[i,bq,cs]$  gives a Tucker approximation of $\F[i,j,k]$ with the same mode-$1$ rank. 
Therefore, we can approximate dominant mode-$1$ subspace of $\F$ by the one of $\U'.$  
The accuracy of resulted approximation is estimated by the following theorem.

\begin{theorem}\label{thm1}
For tensor $\F[i,j,k]$ given by~\eqref{eq:split} it holds
\begin{equation}\nonumber
  \|\F\|_F        \leq \|\U'\|_F          \|V\|_2 \|W\|_2, \quad 
  \|\F\|_\tn      \leq \|\U'\|_\tn        \|V\|_2 \|W\|_2,
\end{equation}
and for mode-$1$ unfoldings $F=F[i,jk]$ and $U'=U'[i,pqcs]$ it holds
\begin{equation}\nonumber
  \|F\|_2   \leq \|U'\|_2   \|V\|_2 \|W\|_2.
\end{equation}
\begin{proof}
The first and last parts follow directly from matrix inequalities
\begin{equation}\nonumber
 \begin{split}
  \|\F[i,j,k]\|_F & = \|F[i, jk]\|_F \leq \|U'[i,bqcs]\|_F \| W[k,cs] \o V[j,bq] \|_2 = \|\U'\|_F \|V\|_2 \|W\|_2, \\
  \|F[i,jk]\|_2   &                  \leq \|U'[i,bqcs]\|_2 \| W[k,cs] \o V[j,bq] \|_2 = \|U'[i,bqcs]\|_2 \|V\|_2 \|W\|_2.
 \end{split}
\end{equation}
Second part reads
\begin{equation}\nonumber
 \begin{split}
 \|\F[i,j,k]\|_\tn & = \max_{\|\u\|=\|\v\|=\|\w\|=1} \f{\U'[i,pq,cs] \x_2 V[j,pq] \x_3 W[k,cs]}{\u[i] \o \v[j] \o \w[k]} = \\
               {}  & = \max_{\|\u\|=\|\v\|=\|\w\|=1} (\v^\t V)[bq]  (\U' \x_1 \u^\t)[bq,cs] (W^\t \w)[cs] \leq \\
               {}  & \leq \max_{\|\u\|=\|\v\|=\|\w\|=1} \|V^\t \v\| \|\U' \x_1 \u^\t\|_2 \|W^\t \w\| = \\
               {}  & = \left(\max_{\|\u\|=1} \|\U' \x_1 \u^\t\|_2\right) \left(\max_{\|\v\|=1} \|V^\t \v\|\right) \left(\max_{\|\w\|=1} \|W^\t \w\|\right) 
                     = \|\U'\|_\tn \|V\|_2 \|W\|_2.
 \end{split}
\end{equation}
\end{proof}
\end{theorem}
\begin{corollary}
For certain perturbation $\Delta\U'$ of tensor $\U',$ the corresponding perturbation  $\Delta\F$ can be estimated as follows
\begin{equation}\label{eq:acc}
 \begin{split}
  \frac{\| \Delta\F\|_F  }{\|\F\|_F}   \leq c_F   \frac{\| \Delta\U'\|_F  }{\|\U'\|_F},  & \qquad c_F  =\frac{\|\U'\|_F   \|V\|_2 \|W\|_2 }{\|\F\|_F}; \\
  \frac{\| \Delta F\|_2  }{\| F\|_2}   \leq c_2   \frac{\| \Delta U'\|_2  }{\| U'\|_2},  & \qquad c_2  =\frac{\| U'\|_2   \|V\|_2 \|W\|_2 }{\| F\|_2}; \\
  \frac{\| \Delta\F\|_\tn}{\|\F\|_\tn} \leq c_\tn \frac{\| \Delta\U'\|_\tn}{\|\U'\|_\tn},& \qquad c_\tn=\frac{\|\U'\|_\tn \|V\|_2 \|W\|_2 }{\|\F\|_\tn}.
 \end{split}
\end{equation}
\end{corollary}
\begin{remark} 
For any tensor, $\|\A[i,j,k]\|_\tn \leq \|A[i,jk]\|_2 \leq \|A[i,j,k]\|_F.$
\end{remark}

To find a dominant mode-$1$ subspace of $\U'[i,bq,cs],$ we can use proper decomposition of Gram matrix of the unfolding $U'[i,bqcs],$ that reads
\begin{equation}\label{eq:gram}
 \begin{split}
 A[i,i'] = (U' {U'}^\t)[i,i'] = U[i,ap] \left( \hat G[p,p'] \o \hat H[a,a']  \right) U[a'p',i'], \\
 \hat G[p,p'] = G[p,qs] G[qs,p'], \quad \hat H[a,a'] = H[a,bc] H[bc,a'].
 \end{split}
\end{equation}
Tensor $\U'$ has a simpler structure than $\F,$ and computation of the Gram matrix~\eqref{eq:gram} is faster than~\eqref{eq:ff}. 
However, evaluation of $A[i,i']$ as full $n_1 \times n_1$ array leads to $\O(n^2r^3)$ complexity. 
Looking for the methods with linear in mode size complexity, we are to use the cross approximation algorithms.

\section{Cross approximation of Gram matrices}
Truncated singular/proper decomposition is used in cases where  low-rank approximation is required.
This problem can be solved by faster methods, for example, those based on \emph{cross approximation} 
$A[i,j] \approx \tilde A[i,j] = U[i,J](A[I,J])^{-1}A[I,j],$ where $I$ and $J$ contain indices of certain rows and columns of $A.$ 
This approximation is exact on the \emph{cross} formed by rows $I$ and columns $J,$ but the overall accuracy depends heavily on the properties of $A[I,J].$ 
In~\cite{gt-psa-1995,gtz-psa-1997,gt-maxvol-2001} it is shown that a good choice for $A[I,J]$  is  \emph{maximum volume} (modulus of determinant) submatrix. 
Search of this submatrix in general case is NP-hard problem~\cite{bartholdi-1982}, and alternatives should be used, see~\cite{tee-cross-2000,gostz-maxvol-2010}. 
If the supported cross is iteratively widened at each step by one row and column that intersect on element where residual is maximum in modulus, cross approximation method is equivalent to Gau{\ss}ian decomposition with complete pivoting. 
For Gram matrix the pivot is always on the diagonal and cross approximation is equivalent to unfinished Cholesky decomposition. 
The resulted algorithm exploiting structure of~\eqref{eq:gram} is summarized in Alg.~\ref{alg}.

\begin{algorithm}[ht]
\caption{Cross approximation for Gram matrix~\eqref{eq:gram}} \label{alg}
\begin{algorithmic}[1]
\REQUIRE Structured tensor $\F = \Kron(\G, \H) \x_1 U \x_2 V \x_3 W,$ see~\eqref{eq:f}
\ENSURE  Approximation $\tilde A = X \Lambda X^\t$ for Gram matrix~\eqref{eq:gram}, such that $\|A - \tilde A\|_F \lesssim \eps\|A\|_F$
\item[\textbf{Initialization:}]   $p=0, \quad \tilde A = 0$ 
\STATE  $\hat G[p,p'] = G[p,qs] G[qs,p'], \quad \hat H[a,a'] = H[a,bc] H[bc,a']$    \hfill {$\O(r^{4})$}
\FOR[Compute diagonal of matrix]{$i=1,\ldots,n$}
      \STATE $U_i[a,p] = U[i,ap], \quad d(i) = \f{U_i[a,p] G[p,p']}{H[a,a']U_i[a',p']}$   \hfill {$\O(r^3)$}
\ENDFOR
\STATE $\nrm := \|\d\|_1$
\REPEAT
\STATE $i_\star := \arg \max_i |d(i)| $ \COMMENT{Find new pivot}                          \hfill {$\O(n)$}
\STATE $\a(:,i_\star) := U[:,ap] (H[a,a'] U_{i_\star}[a',p'] G[p',p])[ap]$                \hfill {$\O(nr^2 + r^3)$}
\STATE $\tilde\a(:,i_\star) = X \Lambda (\xx(i_\star, :))^\t$                             \hfill {$\O(np)$}
\STATE $\xx_\star := (\a-\tilde\a) / \sqrt{(a - \tilde a)(i_\star,i_\star)}$              \hfill {$\O(n)$}
\STATE $\d[i] := \d[i] - |\xx_\star[i]|^2$ \COMMENT{Update diagonal of residual}          \hfill {$\O(n)$}
\STATE $\xx_\star =: [X \: \xx'] \b$  \COMMENT{Orthogonalize $\xx_\star$ to $\Span X$}    \hfill $\O(np)$ 
\STATE $\Lambda + \b^\t\b =: V D V^\t$ \COMMENT{Re-diagonalize decomposition}             \hfill $\O(p^3)$
\STATE $X := [X \: \xx'] V, \quad \Lambda := D, \quad \tilde A = X \Lambda X^\t, \quad \err := \|\d\|_1$          \hfill $\O(np^2)$
\UNTIL{$\err \leq \eps\nrm$ \textbf{or} $r = \rmax$ }
\end{algorithmic}
\end{algorithm}

It is easy to see that evaluation of $\hat G$ and $\hat H,$ i.e. Gram matrices of the unfoldings $G[p,qs]$ and $H[a,bc],$ requires $\O(r^4)$ operations in three-dimensional case and $\O(r^{d+1})$ in $d$-dimension case. With precomputed $\hat G$ and $\hat H$ every element $a(i,i')$ is computed in $\O(r^3)$ operations and a column $\a(:,i')$ is computed in $\O(nr^2 + r^3)$ operations for three and $d$-dimensional case.  
For the rediagonalization of $\Lambda + \b^\t\b$ matrix we can use algorithm proposed by Demmel (see~\cite{demmel}, Alg.~5.3) that is implemented by the LAPACK procedure \texttt{slaed3} and has complexity~$\O(p^3).$
We conclude that approximation of rank-$r$ dominant mode subspace of Gram matrix~\eqref{eq:gram} in  $d$-dimensional case requires $\O(nr^3 + r^{d+1})$ operations.

The relation between accuracy of cross approximation of Gram matrices and corresponding low-rank approximation of initial matrices is given by the following theorem.
\begin{theorem}\label{thm2}
Consider a matrix $U = \left[\begin{array}{cc} U_1 & U_2 \end{array}\right].$ 
If the corresponding Gram matrix 
$$
A = U^\t U = \left[ \begin{array}{cc} A_{11} & A_{12} \\ A_{21} & A_{22}  \end{array}   \right]
$$ 
allows the cross approximation 
\begin{equation}\nonumber
\left\| A - \left[ \begin{array}{c} A_{11} \\ A_{21} \end{array} \right]  A_{11}^{-1} \left[\begin{array}{cc} A_{11} & A_{12} \end{array} \right]  \right\|_2 \leq \eps \|A\|_2,
\end{equation}
then there exists a matrix $B$ such that
\begin{equation}\label{eq:cu}
\|U - U_1 B^\t \|_2 \leq \sqrt{\eps} \|U\|_2.
\end{equation}
\begin{proof}
Consider $V = - U_1 A_{11}^{-1} A_{12} + U_2$ and write
$$
V^\t V = A_{21} A_{11}^{-1} A_{11} A_{11}^{-1} A_{12} - A_{21} A_{11}^{-1} A_{12} - A_{21} A_{11}^{-1} A_{12} + A_{22} = A_{22} - A_{21} A_{11}^{-1} A_{12}.
$$
Cross approximation is exact on the selected rows and columns
\begin{equation}\label{eq:res}
A - \left[ \begin{array}{c} A_{11} \\ A_{21} \end{array} \right]  A_{11}^{-1} \left[\begin{array}{cc} A_{11} & A_{12} \end{array} \right] =
      \left[ \begin{array}{cc} 0 & 0 \\ 0 & A_{22} - A_{21} A_{11}^{-1} A_{12} \end{array} \right] = 
      \left[ \begin{array}{cc} 0 & 0 \\ 0 & V^\t V \end{array} \right],
\end{equation}
and it follows that $\|V^\t V\|_2 \leq \eps \|U^\t U\|_2$ and $\|V\|_2 \leq \sqrt{\eps} \|U\|_2.$
We conclude that $B^\t = \left[\begin{array}{cc}I & A_{11}^{-1} A_{12}\end{array}\right]$ provides~\eqref{eq:cu}.
\end{proof}
\end{theorem}
\begin{remark}
For $U$ with $U_1^\t U_1 = I, \: U_2^\t U_2 = \eps I, \: U_1^\t U_2 = 0,$ inequality~\eqref{eq:cu} is sharp.
\end{remark}
\begin{remark}
For fixed $U_1,$ matrix  $B^\t = \left[\begin{array}{cc}I & A_{11}^{-1} A_{12}\end{array}\right] = (U_1^\t U_1)^{-1} U_1^\t U$  provides minimal residual $U - U_1 B^\t$  in Frobenius and spectral norms. See~\cite{gor-cross-2008}, where a nice estimates for accuracy of cross approximation of matrices and tensors are also given.
\end{remark}
\begin{remark}
$\Span B = \Span X$ is the subspace of columns of the Gram matrix that support the cross approximation in Alg.~\ref{alg}.
\end{remark}

Since the spectral norm of the residual is not easy to evaluate, the stopping criteria in a practical algorithm is based on the Frobenius norm.
On each step of Alg.~\ref{alg} vector $\d$ contains the diagonal of residual~\eqref{eq:res} and
$$
\|\d\|_1 =  \sum_{i=1}^n |d(i)| = \sum_{i} (V^\t V)(i,i) = \|V[i,j]\|_F^2.
$$
We can also implement stopping criteria based on eigenvalues stored in $\Lambda.$ To do this, we can split them in `dominant' and `smaller' parts basing on desired tolerance $\eps,$ and stop the process if during several iterations new eigenvalues fall into the smaller part.  
This criteria will approximate spectral norm more precisely, but as we see in numerical experiments, it generally does not differ from the Frobenius-based one.

Obviously, Alg.~\ref{alg} can be applied in the same way to estimate other Tucker factors of~\eqref{eq:f}. 
Due to roundoff errors, accuracy $\eps$ of Alg~\ref{alg} is limited by machine precision $\tol,$ and for $\eps=\tol$, accuracy of~\eqref{eq:appr} can be estimated by Thm.~\ref{thm2} as
$$
\|\F - \tilde\F\|_2 \leq \sqrt{\tol} \sqrt{c_2^2(U) + c_2^2(V) + c_2^2(W)} \|\F\|_2,
$$
where $c_2(U)$ is defined in~\eqref{eq:acc} and similar definition applies to $V, W.$

\section{Numerical examples}
Multidimensional data often appear in modern modelling programs. 
For example, in chemical packages, e.g. PC GAMESS, MOLPRO, the \emph{electron density function} is given in canonical form~\eqref{eq:c} as a sum of tensor product of Gaussians, but with number of terms, that may be too large for practically feasible computations even for moderate molecules. 
In order to make computations efficient,  further approximation (recompression) to the Tucker format can be performed. 
This problem was approached in~\cite{mpi-chem3d-2009} using Tucker-ALS algorithm, in \cite{khor-ml-2009} by Tucker-ALS with initial guess obtained from the coarser grids, in~\cite{fkst-chem-2008} by Cross3D algorithm, in~\cite{ost-chem-2009} by individual cross approximation of canonical factors, in~\cite{sav-rr-2009} by cross approximation of Gram matrices of unfoldings and in~\cite{gos-kryl-2010} by algorithms based on Wedderburn rank reduction.

As an example, we apply the discussed algorithm for Hadamard multiplication of electron density given in Tucker format to themselves. 
This operation can be a building block for algorithm that computes pointwise cubic root of density, that is used in the Kohn-Sham model. 
A good initial guess for such methods can be evaluated by mimic algorithm~\cite{ost-chem-2009}. 

The results of experiments are collected in Tab.~\ref{tab}. 
They were performed on Intel Xeon Quad-Core E5504 CPU running at $2.00$~GHz using Intel Fortran compiler version 11.1 and BLAS/LAPACK routines provided by MKL library.
For each molecule, we show time in seconds $T(\mbox{Alg.~\ref{alg}})$ for evaluation of three dominant subspaces 
$X[i,\alpha], Y[j,\beta]$ and $Z[k,\gamma]$  by Alg.~\ref{alg} with accuracy of approximation of Gram matrices set to $\eps=10^{-12}.$ 
Then we compute best core by convolution
\begin{equation}\nonumber
 \T[\alpha,\beta,\gamma] = \F[i,j,k] \x_1 X[\alpha,i] \x_2 Y[\beta,j] \x_3 Z[\gamma,k].
\end{equation}
and check relative accuracy $\eps(\mbox{Alg.~\ref{alg}})$ of approximation~\eqref{eq:appr} in Frobenius norm. 
The direct computation of all elements of residual requires a lot of computational time and the accuracy $\|\F-\tilde\F\|_F$ was verified by comparing the result with Tucker approximation computed by Cross3D algorithm~\cite{ost-tucker-2008} with accuracy set to $\eps=10^{-12}.$ The Cross3D algorithm was verified in~\cite{ost-tucker-2008,fkst-chem-2008} by exhaustive check on parallel memory platforms, and can be considered as reliable answer. The residual between two Tucker formats is computed as proposed in~\cite{ost-latensor-2009}.

Then we compute approximation of the same accuracy~$\eps(\mbox{Alg.~\ref{alg}})$ by Cross3D~\cite{ost-tucker-2008} and WsvdR~\cite{gos-kryl-2010} algorithms and show the corresponding timings as $T(\mbox{c3d})$ and $T(\mbox{wsvdr})$. 
We also show time $T(\mbox{tals})$ for one iteration of Tucker-ALS~\cite{tuckerals-1980,lathauwer-rank1-2000} with ranks fixed equal to the ranks of bases $X, Y, Z,$ returned by Alg.~\ref{alg}.  
Then we apply one iteration of rank-revealing Tucker-ALS~\cite{ost-chem-2009} with accuracy parameter set to $\eps=10^{-12}$ using bases $X, Y, Z$ as initial guess, and show the accuracy of improved approximation by~$\eps(\mbox{tals}).$

\begin{table}[t]
\caption{Hadamard square of electron density, $n_1=n_2=n_3=5121$} \label{tab}
 \begin{center}
  \begin{tabular}{cc|cc|ccc|c}
   molecule & $r_1, r_2, r_3$ & $T(\mbox{Alg.~\ref{alg}})$ & $\eps(\mbox{Alg.~\ref{alg}})$ & $T(\mbox{c3d})$ & $T(\mbox{wsvdr})$ & $T(\mbox{tals})$ & $\eps(\mbox{tals})$ \\ \hline
   methane & $(74,74,74)$    & $4.0$ & $3\eee{-7}$ & $78.6$ & $12.4$ & $37$  & $7\eee{-13}$ \\
   ethane  & $(67,94,83)$    & $5.3$ & $6\eee{-7}$ & $76.8$ & $15.1$ & $42$  & $8\eee{-13}$ \\
   ethanol & $(128,127,134)$ & $20$  & $5\eee{-7}$ & $1050$ & $210$  & $473$ & $9\eee{-13}$  \\
   glycine & $(62,176,186)$  & $38$  & $8\eee{-7}$ & $1260$ & $237$  & $442$ & $9\eee{-13}$ \\
  \end{tabular}
 \end{center}
\end{table}

We conclude that proposed algorithm is faster that other methods for this purpose and return approximation of dominant subspaces that allows to construct approximation with accuracy about square root of machine precision. Using the subspaces, computed by Alg.~\ref{alg} as initial guess, rank revealing Tucker-ALS converges to almost machine precision in one iteration.

\section*{Acknowledgements}
This work was supported by RFBR grants 08-01-00115, 09-01-12058, 10-01-00757, RFBR/DFG grant 09-01-91332, Russian Federation Gov. contract $\Pi940$ and  Priority research program of Dep. Math. RAS.
The first author was supported by RFBR travel grant 10-01-09201 to present the results of this paper on ICSMT (Hong Kong, January 2010).
Part of this work was done during the stay of the first author in Max-Plank Institute for Mathematics in Sciences in Leipzig (Germany). 
Authors are grateful to Heinz-J\"urgen Flad and Rao Chinnamsettey for providing input data for the electron density functions.


\end{document}